\documentclass[10pt, leqno]{amsart}
\usepackage{amssymb,amsfonts,amscd,amsbsy}
\usepackage[mathscr]{eucal}

\theoremstyle{plain}
\newtheorem{thm}{Theorem}[section]

\newtheorem{pro}[thm]{Proposition}
\newtheorem{cor}[thm]{Corollary}
\theoremstyle{definition}

\newtheorem{rem}[thm]{Remark}

\begin{document}

\title[Cyclotomic fields]{Vanishing of eigenspaces and Cyclotomic fields}
\author[Robert Osburn]{Robert Osburn}

\address{Max-Planck Institut f{\"{u}}r Mathematik, Vivatsgasse 7, 53111 Bonn,
Germany}

\email{osburn@mpim-bonn.mpg.de}
 
\subjclass[2000]{Primary: 11R18; Secondary: 11T22}

\begin{abstract}
We use a result of Thaine to give an alternative proof of the fact that, for a
prime $p>3$ congruent to $3$ modulo $4$, the component $e_{(p+1)/2}$ of the
$p$-Sylow subgroup of the ideal class group of $\mathbb Q(\zeta_{p})$ is
trivial. 
\end{abstract}
 
\maketitle
\section{Introduction}
Let $p>3$ be a prime, $\zeta_{p}$ a $p$th primitive root of $1$, and $\Delta$
the Galois group of $\mathbb Q(\zeta_{p})$ over $\mathbb Q$. Let $q\neq p$ be a
prime and $n$ the order of $q$ modulo $p$. Assume $q \not\equiv 1 \bmod p$ and so 
$n \geq 2$, $p(q-1)|q^n-1$, and $n|p-1$. Set $f=(q^n-1)/p$ and $e=(p-1)/n$. 
Let $Q$ be a prime ideal of $\mathbb Z[\zeta_{p}]$ above $q$ and let $\mathbb
F=\mathbb Z[\zeta_{p}]/Q$. Thus $\mathbb F \cong \mathbb F_{q^n}$, the finite field with
$q^n$ elements. Let $\alpha \in \mathbb Z[\zeta_{p}]$ be a generator of $\mathbb
F^{\times}$ such that

\begin{center}
$\alpha^{f} \equiv \zeta_{p} \bmod Q$.
\end{center}

Now let $A$ be the $p$-Sylow subgroup of the ideal class group $\mathbb
Q(\zeta_p)$, $\mathbb Z_{p}$ the ring of $p$-adic integers, $\omega: \Delta \to
\mathbb Z_{p}^{\times}$ is the Teichm{\"u}ller
character defined by 

\begin{center}
$\omega(k) \equiv k \bmod p$,
\end{center} 

\noindent and $e_{r}$, $0 \leq r \leq p-2$, the idempotents
 
\begin{center}
$\displaystyle \frac{1}{p-1} \sum_{\lambda \in \Delta}^{}
\omega^{r}(\lambda) {\lambda}^{-1} \in \mathbb Z_{p}[\Delta]$.
\end{center} 
 
As $A$ is a $\mathbb Z_{p}[\Delta]$-module, we have the decomposition (see
Section 6.3 in \cite{w})

\begin{center}
$\displaystyle A=\bigoplus_{r=0}^{p-2} e_{r}(A)$.
\end{center}

It is well-known that for $r$ even, $2 \leq r \leq p-3$, $e_{r}(A)$ can be
identified with the components of the $p$-part of the ideal class group of
$\mathbb Q(\zeta_{p} + \zeta_{p}^{-1})$. Vandiver's conjecture says that all
even components $e_{r}(A)$ vanish. Via K-theory, Kurihara \cite{k} proved that the ``top''
even eigenspace $e_{p-3}(A)$ always vanishes. Kurihara's proof uses the
surjectivity of the Chern map

\begin{center}
$K_{4}(\mathbb Z) \otimes \mathbb Z/p \to e_{p-3}(A)$.
\end{center}

Soul\'e \cite{s} extended Kurihara's result and showed the following: Let $n>1$
odd. If $\log p> n^{224n^4}$, then $e_{p-n}(A)$ is trivial. Our main result is the
following.

\begin{thm}
If $p>3$ is a prime congruent to $3$ modulo $4$, then $e_{(p+1)/2}(A)$ is trivial.
\end{thm}

One can use the reflection theorem (see Theorem 10.9 in \cite{w}) and the class
number formula for the imaginary quadratic field $\mathbb Q(\sqrt{-p})$ to prove
Theorem 1.1. Precisely, if $e_{(p+1)/2}(A)$ is non-trivial, then $e_{(p-1)/2}(A)$
is non-trivial and so (see Section 2 for a precise description of $v$) $p\mid n-2v$;
but $n-2v<p$, a contradiction. Recently, Thaine \cite{th} investigated properties of
certain numbers (see $d_{i}$ and $a_{k}$ in Section 2) related to Gaussian periods and
showed that they are useful in the study of certain components of the ideal class
group of $\mathbb Q(\zeta_{p})$. In particular, Thaine states: ``We believe that
the theorem [Theorem 1 in \cite{th} or Theorem 2.2 below] can be used to show that, with $l$ odd ($1 \leq l \leq e-1$), some of
the components $e_{p-ln}(A)$ of $A$ are trivial. The idea is to show that if
$e_{p-ln}(A)$ is non-trivial, then {\it all} prime numbers $q$ of order $n$ modulo $p$
must have a certain form; we hope this will contradict some version of Dirichlet's
theorem on primes in arthimetic progressions.''

The purpose of this note is to give an alternative proof of Theorem 1.1 using Thaine's
result. The proof will show that non-trivial eigenspaces lead
to representations of certain integers by binary quadratic forms with
restrictive divisibility properties on the parameters. A density calculation will
show that that this divisibility property doesn't occur for all primes $q$ of order
$(p-1)/2$ modulo $p$ and thus the eigenspace must vanish. It might be of some
interest to see if similar vanishing results can be obtained for other even
indexed eigenspaces $e_{r}(A)$ by considering an appropriate quadratic form (see
(2) of Remark 4.1).

\section{Indices of Cyclotomic Units}

We discuss a result of Thaine on congruences involving indices of cyclotomic
units. Let us consider the components $e_{p-ln}(A)$ for $l$ odd, $1 \leq l \leq
e-1$. For $r$ even, $2 \leq r \leq p-3$, let

\begin{center}
$\displaystyle \beta_{r}= \prod_{i=1}^{p-1} (1-\zeta_{p}^{i})^{i^{p-1-r}}$
\end{center}

\noindent and let $i_{r}(Q)$ be the least nonnegative integer such that

\begin{center}
$\beta_{r} \equiv \alpha^{i_{r}(Q)} \bmod Q$.
\end{center}
 
It is well-known that $e_{r}(A)$ is trivial if and only if $\beta_{r}$ is not
the $p$th power of an element of $\mathbb Z[\zeta_{p}]$ (see Theorem 15.7 and
the discussion preceding Theorem 8.14 in \cite{w}). In particular, we have the
following.

\begin{pro} For $r$ even, $2 \leq r \leq p-3$, if $i_{r}(Q) \not\equiv 0 \bmod
p$, then $e_{r}(A)$ is trivial.
\end{pro}
  
In order to study the indices $i_{r}(Q)$, we need to introduce certain numbers.
Let $g$ be a primitive root modulo $p$. For $k \in \mathbb Z$, we define

\begin{center}
$\displaystyle a_{k}= nq^{v}\sum_{i=0}^{e-1} g^{nki}d_{i}$,
\end{center} 

\noindent where $d_{i}$ is defined (see (14) in \cite{th}) as 

\begin{center}
$\displaystyle d_{i}=\frac{\eta_{g^{i}}- \eta_{0}}{q^v}$
\end{center}

\noindent for $0 \leq i \leq e-1$ where the $\eta_{m}$'s are the Gaussian
periods $(0 \leq m \leq p-1)$

\begin{center}
$\displaystyle \eta_{m}=\sum_{j=0}^{f-1} \zeta_{q}^{T({\alpha}^{m+pj})}$.
\end{center}

\noindent Here $T$ is the trace from $\mathbb F$ to $\mathbb F_q$. Given an
integer $a$, denote by $|a|_p$ the smallest nonnegative residue of $a$ modulo
$p$. Let $q^{v}$ be the largest power of $q$ dividing a certain Gauss sum
$G(\zeta_{p})$ (see p. 314 in \cite{th}). It follows from (12) in \cite{th} that

\begin{center}
$\displaystyle v=\underset{0 \leq k \leq e-1}{\text{min}} \frac{1}{p} \sum_{l=0}^{n-1}
|g^{k+el}|_{p}$.
\end{center}

\noindent Note that $v\geq 1$. The following theorem (see Theorem 1 in
\cite{th}) summarizes some properties of the numbers $a_{k}$ and $d_{i}$ and
will be useful when considering $e_{p-ln}(A)$.
 
\begin{thm} 
(i) We have 

\begin{center}
$\displaystyle e^{2}q^{n-2v}=\Bigg(\sum_{i=0}^{e-1} d_i \Bigg)^2 +
p\Bigg( e\sum_{i=0}^{e-1} d_{i}^2 - \Big(\sum_{i=0}^{e-1} d_i \Big)^2 \Bigg)$.
\end{center}

(ii) The numbers $a_k$ satisfy the following congruences: $a_{0} \equiv -1 \bmod
p$, and for $l$ odd $(1 \leq l \leq e-1)$, 

\begin{center}
$\displaystyle \sum_{m=1}^{l} (-1)^m \binom{ln}{mn}a_{l-m}a_{m} \equiv
-l \cdot i_{p-ln}(Q) \bmod p$.
\end{center}
\end{thm}  

\section{Quadratic Forms}
Let $p>3$ be a prime congruent to $3$ modulo $4$. In this section, we study the
representation of primes and a multiple of a certain power of a prime by the quadratic form $x^2+py^2$. 

\begin{pro}
If $p>3$ is a prime and $p \equiv 3 \bmod 4$, then there exists a prime $q\neq p$ with
$\Bigl(\frac{-p}{q}\Bigr)=1$ and integers $u$ and $v$ such that $p \nmid v$ and
$$u^2+pv^2=q.$$
\end{pro}

\begin{proof}
By Theorem 9.12 in \cite{cox}, there are infinitely many primes $q$ with
$\Bigl(\frac{-p}{q}\Bigr)=1$ such that $x^2+py^2=q$ has an integer solution.
Let $\mathcal{S}_1$ denote the set of primes represented by $x^2+py^2$ and
$\mathcal{S}_2$ denote the set of primes represented by $x^2+p^3y^2$. 
Suppose for every prime $q\neq p$, we have $p \mid v$. Then the quadratic forms
$x^2+py^2$ and $x^2+p^3y^2$ represent the same infinite set of primes and thus
$\mathcal{S}_1$ and $\mathcal{S}_2$ have the same Dirichlet density. 
Let $h$ be the class number of $\mathbb Q(\sqrt{-p})$. By Theorems 7.24 and 9.12 in \cite{cox}, we have the following: for $p \equiv 7 \bmod 8$, 

\begin{center}
$\mathcal{S}_1$ has density $\displaystyle \frac{1}{2h}$
\end{center}

\noindent and 

\begin{center}
$\mathcal{S}_2$ has density $\displaystyle \frac{1}{2ph}$,
\end{center}

\noindent which is a contradiction. For $p \equiv 3 \bmod 8$, 

\begin{center}
$\mathcal{S}_1$ has density $\displaystyle \frac{1}{6h}$
\end{center}

\noindent and 

\begin{center}
$\mathcal{S}_2$ has density $\displaystyle \frac{1}{6ph}$, 
\end{center}

\noindent a contradiction. Therefore, there exists a prime $q$ such that $p\nmid v$.
\end{proof}

We now need the following result from \cite{bh} (see Theorem 2, page 224).

\begin{thm} Let $p$ and $q$ be distinct odd primes and assume that $u^2+pv^2=q$ for integers $u$ and
$v$. If $s$ is an odd positive integer, let $x(U,V)$ and $y(U,V)$ be the
polynomials in $\mathbb Z[U,V]$ such that

\begin{center}
$(U+V\sqrt{-p})^{s} = x(U,V) + y(U,V)\sqrt{-p}$.
\end{center}

\noindent Then 

\begin{equation}
\displaystyle y(u,v)=v\sum_{j=0}^{(s-1)/2} \binom{s}{2j}(-pv^2)^{(s-2j-1)/2}(u^2)^j
\end{equation}

\noindent and ($x(u,v)$, $y(u,v)$) is a solution to $x^2+py^2=q^s$.
\end{thm}

\begin{cor} Let $p>3$ be a prime with $p \equiv 3 \bmod 4$ and $h$ be the class
number of $\mathbb Q(\sqrt{-p})$. Then there a prime $q \neq p$ with
$\Bigl(\frac{-p}{q}\Bigr)=1$ such that if $(C,D)$ is a solution of the equation

\begin{center}
$x^2+py^2=4q^h$,
\end{center}

\noindent then $p\nmid D$. 
\end{cor}  
  
\begin{proof}
By Proposition 3.1, there exists a prime $q\neq p$ such that for the integer solution
$(u,v)$ to the equation $x^2+py^2=q$, we have $p \nmid v$. We also have that $p
\nmid u$ since otherwise $p$ would divide $q$. Now take $s=h$ in Theorem 3.2. By
(1) and the fact that $h$ is odd (see Corollary 18.4 in \cite{cnp}), we see that
for the solution ($x(u,v)$, $y(u,v)$) to the equation 

\begin{equation}
x^2+py^2=q^h,
\end{equation}

\noindent $p \nmid y(u,v)$. Multiplying (2) by 4, we have $p \nmid D$.
\end{proof}

\section{Proof of Theorem 1.1}

\begin{proof}
Let $p>3$ be a prime with $p \equiv 3 \bmod 4$, $q \neq p$ be a prime of order
$n=(p-1)/2$ modulo $p$, and $g$ a primitive root modulo $p$. Let us note that for a fixed prime $p>3$ with $p\equiv 3 \bmod 4$, a prime $q$
of order $(p-1)/2$ modulo $p$ always exists. To see this, consider the group $G=(\mathbb Z/p)^{*}$. This
is a cyclic group of order $p-1$. There exists a unique element $\alpha_{p}$ in $G$ of order $(p-1)/2$. By
Dirichlet's theorem on primes in arithmetic progressions, there are infinitely
many primes $q$ such that $q \equiv \alpha_{p} \bmod p$. 

Now if $q$ has order $(p-1)/2$ mod $p$, then $\Bigl(\frac{q}{p}\Bigr)=1$. We have that
$\Bigl(\frac{-p}{q}\Bigr)=\Bigl(\frac{-1}{q}\Bigr)\Bigl(\frac{p}{q}\Bigr)$. If $q\equiv 1 \bmod 4$, then $\Bigl(\frac{-1}{q}\Bigr)=1$ and
$\Bigl(\frac{p}{q}\Bigr)=\Bigl(\frac{q}{p}\Bigr)=1$. Thus
$\Bigl(\frac{-p}{q}\Bigr)=1$. If $q\equiv 3 \bmod 4$, then $\Bigl(\frac{-1}{q}\Bigr)=-1$ and
$\Bigl(\frac{p}{q}\Bigr)=-1 \iff \Bigl(\frac{q}{p}\Bigr)=1$. Thus
$\Bigl(\frac{-p}{q}\Bigr)=1$. Thus if $q$ is a prime of order $(p-1)/2$ mod $p$, then
$\Bigl(\frac{-p}{q}\Bigr)=1$. Suppose $e_{(p+1)/2}(A)$ is nontrivial. Then by
Proposition 2.1, $i_{(p+1)/2}(Q) \equiv 0 \bmod p$. By Theorem 2.2, $(ii)$, we 
have that 

\begin{center}
$a_1 \equiv 0 \bmod p$ 
\end{center}

\noindent and so by the definition of $a_1$, 

\begin{center}
$d_1 \equiv d_0 \bmod p$. 
\end{center}

\noindent By Theorem 2.2, $(i)$,

$$
\begin{aligned}
4q^{(p-1)/2 - 2v} &= (d_0 + d_1)^2 + p(2(d_{0}^{2} + d_{1}^{2})-(d_0 - d_1)^{2})\\
&=(d_0 + d_1)^2 + p(d_0 - d_1)^2,
\end{aligned}
$$

\noindent where $p\mid (d_0 - d_1)$. We now claim that $(p-1)/2-2v=h$ where $h$ is the class number of $\mathbb
Q(\sqrt{-p})$. To see this, note that as $n=(p-1)/2$, then $e=2$ and $l=1$. Thus 

\begin{center}
$\displaystyle v=\text{min} \Bigg\{\frac{1}{p} \sum_{l=0}^{\frac{p-3}{2}}
|g^{2l}|_{p}, \frac{1}{p} \sum_{l=0}^{\frac{p-3}{2}}
|g^{1+2l}|_{p} \Bigg\}$. 
\end{center}

\noindent As $g$ is a primitive root modulo $p$, the first sum, say $R$, appearing in $v$ is

\begin{center}
$\displaystyle \frac{1}{p} \sum_{\substack{s=0 \\\bigl(\frac{s}{p}\bigr)=1}}^{\frac{p-1}{2}}
s$. 
\end{center}
 
\noindent Similarly, the second sum, say $V$, appearing in $v$ is

\begin{center}
$\displaystyle \frac{1}{p} \sum_{\substack{t=0
\\\bigl(\frac{t}{p}\bigr)=-1}}^{\frac{p-1}{2}} t$.
\end{center}

\noindent By \cite{dir}, $V+R=(p-1)/2$ and $h=V-R$. This implies that $v=R$. Thus $R=(p-1)/4-h/2$ and so

\begin{center}
$(p-1)/2 - 2v=(p-1)/2 - 2R = h$.
\end{center}

\noindent So for every prime $q$ of order $(p-1)/2$ mod $p$ (and thus $\Bigl(\frac{-p}{q}\Bigr)=1$), if $a:=(d_0 + d_1)$ and $b:=(d_0 - d_1)$ are the integers such that

\begin{center}
$4q^{h}=a^2 + pb^2$,
\end{center}

\noindent then we have that $p \mid b$. By Corollary 3.3, this is a contradiction. Thus $e_{(p+1)/2}(A)$ is trivial.
\end{proof}

\begin{rem} (1) The quadratic form which appears in the proof of Theorem 1.1 has been studied 
by Stickelberger \cite{stick}. His elegant result is as follows:
Let $-k$ be a negative fundamental discriminant so that $\mathbb Q(\sqrt{-k})$
is an imaginary quadratic field of discriminant $-k$ and class number $h$.
Assume $k\neq 3$, $4$, or $8$. Let $q$ be a prime such that $\Bigl(\frac{-k}{q}\Bigr)=1$. 
Then there are integers $C$ and $D$, unique up to sign, for which 

\begin{center}
$4q^{h}=C^2 + kD^2$
\end{center}

\noindent where $k \nmid C$. Moreover, for prime $k \geq 7$, $C\equiv 2(-q)^{-R}
\bmod k$ where $R$ is as above.

(2) It might be possible to prove similar vanishing results for $e_{(3p+1)/4}(A)$ where $p \equiv 5 \bmod 8$
and for $e_{(5p+1)/6}(A)$ where $p \equiv 7 \bmod 12$. The difficulty is that the
resulting quadratic forms are more involved. Namely, the first case corresponds
to $e=4$ and so by Theorem 2.2, $(i)$, we have

\begin{center}
$\displaystyle 16q^{(p-1)/4 - 2v} = \Bigg(\sum_{i=0}^{3} d_i \Bigg)^2 + 
p\Bigg(4\sum_{i=0}^{3} d_{i}^2 - \Big(\sum_{i=0}^{3} d_i \Big)^2 \Bigg).$
\end{center}

The second case corresponds to $e=6$ and thus by Theorem 2.2, $(i)$, we have

\begin{center}
$\displaystyle 36q^{(p-1)/6 - 2v} = \Bigg(\sum_{i=0}^{5} d_i \Bigg)^2 + p\Bigg( 6\sum_{i=0}^{5} d_{i}^2 -  
\Big(\sum_{i=0}^{5} d_i \Big)^2 \Bigg).$
\end{center}
\end{rem}

\section*{Acknowledgments}

The author would like to thank Francisco Thaine and Pieter Moree for their
comments and encouragement. Thanks also to the referee for suggestions which
improved the exposition of the paper. Finally, the author would
like to thank the Max-Planck-Institut f{\"u}r Mathematik for their
hospitality and support during the preparation of this paper.

\end{document}